# FALSE DISCOVERY AND FALSE NONDISCOVERY RATES IN SINGLE-STEP MULTIPLE TESTING PROCEDURES[1]

By Sanat K. Sarkar

*Temple University*

Results on the false discovery rate (FDR) and the false nondiscovery rate (FNR) are developed for single-step multiple testing procedures. In addition to verifying desirable properties of FDR and FNR as measures of error rates, these results extend previously known results, providing further insights, particularly under dependence, into the notions of FDR and FNR and related measures. First, considering fixed configurations of true and false null hypotheses, inequalities are obtained to explain how an FDR- or FNR-controlling single-step procedure, such as a Bonferroni or Šidák procedure, can potentially be improved. Two families of procedures are then constructed, one that modifies the FDR-controlling and the other that modifies the FNR-controlling Šidák procedure. These are proved to control FDR or FNR under independence less conservatively than the corresponding families that modify the FDR- or FNR-controlling Bonferroni procedure. Results of numerical investigations of the performance of the modified Šidák FDR procedure over its competitors are presented. Second, considering a mixture model where different configurations of true and false null hypotheses are assumed to have certain probabilities, results are also derived that extend some of Storey's work to the dependence case.

**1. Introduction.** The false discovery rate (FDR) and related measures have been receiving considerable attention due to their relevance as measures of the overall error rate in multiple testing problems that arise in many scientific investigations, particularly in the context of DNA microarray analysis. Consider Table 1, which summarizes the outcomes in multiple testing of $n$ null hypotheses $H_1, \ldots, H_n$. Let $Q = V/R$ if $R > 0$ and $= 0$ if

Received November 2003; revised March 2005.

[1]Supported by NSF Grant DMS-03-06366 and a 2003 Summer Research Fellowship awarded by Temple University.

TABLE 1
*The outcomes in testing n null hypotheses*

|            | Rejected | Accepted | Total |
|------------|----------|----------|-------|
| True null  | V        | U        | $n_0$ |
| False null | S        | T        | $n_1$ |
| Total      | R        | A        | $n$   |

$R = 0$, that is, the proportion of false positives (Type I errors) among the rejected null hypotheses. Genovese and Wasserman [9] called this the false discovery proportion (FDP). The FDR is defined by $E(Q)$. It was first introduced in multiple testing by Benjamini and Hochberg [1], who provided a step-up procedure that controls the FDR with independent test statistics. Later, Benjamini and Liu [4] offered a step-down FDR procedure under independence. The FDR-controlling property of the Benjamini–Hochberg (BH) procedure was extended by Benjamini and Yekutieli [5] to some positively dependent multivariate distributions. Sarkar [14] proved that the critical values of the BH procedure can be used in a more general stepwise procedure to provide control of the FDR not only under independence, but also when the test statistics have the same type of positive dependence property as considered by Benjamini and Yekutieli [5]. In addition, he established the FDR-controlling property of the Benjamini–Liu step-down procedure for some positively dependent test statistics. Genovese and Wasserman [8, 9] investigated some operating characteristics of the BH procedure asymptotically under independence and further extended the theory of FDR by taking a stochastic process approach.

A slightly different concept of FDR, called the positive false discovery rate (pFDR), was considered by Storey [17]. It is defined as the conditional FDR given at least one rejection, that is, $\mathrm{pFDR} = E(V/R|R > 0)$, and it has the interpretation of a Bayesian Type I error rate under a mixture model involving i.i.d. $p$-values when a single-step multiple testing procedure is used; see also [18]. Storey [17] provided estimates of FDR and pFDR under the above mixture model for a single-step procedure that are related to the empirical Bayes FDR of Efron, Tibshirani, Storey and Tusher [7]; see also [6]. A new family of FDR procedures based on estimates of FDR was suggested by Storey [17] and Storey, Taylor and Siegmund [19].

An analog of FDR in terms of false negatives (Type II errors) was introduced by Genovese and Wasserman [8] and Sarkar [15]. It is the FNR, called false nondiscovery rate by Genovese and Wasserman [8] and the false negatives rate by Sarkar [15]. It is defined by $E(N)$, where $N = T/A$ if $A > 0$ and $= 0$ if $A = 0$ is the proportion of false negatives among the accepted null hypotheses or the false nondiscovery proportion (FNP) [9]. Storey [18] defined the pFNR (positive false nondiscovery rate), the conditional expectation $E(T/A|A > 0)$, as an analog of his pFDR. While Genovese and Wasserman [8] considered new methods that incorporate both FDR and FNR, Storey [18] established a connection between multiple testing and classification theory in terms of a combination of pFDR and pFNR. Sarkar [15] proved that the FNR can be controlled by a step-down analog of the BH procedure. He also introduced a concept of unbiasedness of an FDR- or FNR-controlling multiple testing procedure and established this property for a generalized stepwise procedure under independence.

In this article we mainly concentrate on single-step multiple testing procedures, and we develop new results on FDR and FNR with dependent test statistics both under a model where the configuration of true and false null hypotheses is assumed fixed, yet unknown, and under the so-called mixture model where different configurations of true and false null hypotheses are assumed to have certain probabilities. The intent of these results is to verify some desirable properties of FDR and FNR and to extend some previously known results, thereby providing further insights into the notions of FDR and FNR and related measures, particularly under dependence.

Suppose that $\mathbf{X} = (X_1, \ldots, X_n)$ has a joint distribution indexed by the set of parameters $\boldsymbol{\theta} = (\theta_1, \ldots, \theta_n)$. Let $H_i : \theta_i \leq \theta_{i0}$ be tested against $K_i : \theta_i > \theta_{i0}$, for some given $\theta_{i0}$, $i = 1, \ldots, n$. Let $\{H_i : i \in J_0\}$ and $\{H_i : i \in J_1\}$ be the sets of true and false null hypotheses, respectively. It will be assumed that $J_0$ is nonempty. Consider a single-step procedure that rejects $H_i$ in favor of $K_i$ if $X_i \geq t$ for some fixed $t$. Two of our main results with fixed $J_0$ and $J_1$ (Theorems 1 and 3) are that if $\mathbf{X}$ is stochastically increasing in each $\theta_i$, which is typically the case in many multiple testing problems, then the maximum values of FDR and FNR of a single-step procedure are $(n_0/n)P\{R > 0\}$ and $(n_1/n)P\{A > 0\}$, respectively, where the probabilities are evaluated at $\theta_0 = (\theta_{10}, \ldots, \theta_{n0})$ and $\mathbf{X}$ is assumed exchangeable under these null hypothesis values. In addition to representing more precise versions of the results that state that Šidák and Bonferroni single-step procedures control FDR or FNR, these theorems show how these procedures can potentially be improved in terms of having better control of FDR or FNR borrowing information about $n_0$ or $n_1$ from the data in the spirit of Benjamini and Hochberg [2], Benjamini, Krieger and Yekutieli [3], Storey [17] and Storey, Taylor and Siegmund [19]. Storey, Taylor and Siegmund [19] provided procedures for modifying the BH procedure using a family of estimates of $n_0$ and proved that they control FDR under independence. We obtain new families of procedures: one to modify the FDR-controlling and the other to modify the FNR-controlling Šidák procedure. Considering independent test statistics, we prove that they control FDR or FNR. The modified Šidák FDR procedures are less conservative under independence than the corresponding family that modifies the Bonferroni procedure obtained by using the estimates of $n_0$ considered in [19]. An analogous result is true for modified Šidák FNR procedures. Our method of modifying the Šidák FDR and the Šidák FNR procedures relies directly on two new results, Theorems 2 and 4, which extend inequalities given by Theorems 1 and 3, respectively, under independence from a single-step to a two-step procedure.

Next, we derive certain results that extend Storey's [17, 18] work to the dependent case. Storey obtained expressions for the FDR and FNR of a single-step procedure under a mixture model where, given any configuration of true and false null hypotheses, the $X_i$'s are assumed to be independent, providing useful Bayesian interpretations to his notions of pFDR and pFNR. More specifically, he proved: pFDR $= P\{H_1$ is true$|X_1 \geq t\}$ and pFNR $= P\{H_1$ is false$|X_1 < t\}$, irrespective of the number of tests. Assuming a more general mixture model in which the $X_i$'s are assumed to be dependent with a location family of distributions and to have a certain type of positive dependence structure, we prove in Theorems 5 and 6, respectively, that pFDR $\leq \max_{1 \leq i \leq n} P\{H_i$ is true$|X_i \geq t\}$ and pFNR $\leq \max_{1 \leq i \leq n} P\{H_i$ is false$|X_i < t\}$, with the equalities holding under independence. An important implication of the first inequality is that Storey's [17] $q$-value for a single-step multiple test under certain commonly encountered types of dependence is more conservative, as one would desire, than that under independence.

The paper is organized as follows. In Section 2 we formally define the stochastic increasing property we need for $\mathbf{X}$ to obtain the maximum values of FDR and FNR for fixed $J_0$ and $J_1$. Section 3 reports the results related to FDR for fixed $J_0$ and $J_1$, and some numerical results that show the performance of the modified Šidák procedure in controlling FDR compared to the modified Bonferroni and the original Bonferroni and Šidák procedures. Similar results related to FNR are presented in Section 4, of course without showing any additional numerical evidence. Section 5 numerically compares the Bonferroni and Šidák procedures with their modified versions in terms of a concept of power involving both FDR and FNR. Section 6 presents the results on FDR and FNR under the aforementioned mixture model with dependent $\mathbf{X}$. Proofs are given in Section 7. The paper concludes with some final remarks in Section 8.

**2. Stochastically increasing family of distributions.** This section defines a type of stochastic increasing property of a family of distributions that will be required to establish our results on FDR and FNR. Whenever an increasing or decreasing condition or property in terms of $\mathbf{X}$ or $\theta$ is mentioned, it is to be understood as being coordinatewise.

DEFINITION 1. An $n$-dimensional random vector $\mathbf{X} = (X_1, \ldots, X_n)$ or the corresponding family of distributions $\{P_\theta\}$, where $\theta = (\theta_1, \ldots, \theta_n)$, is said to be stochastically increasing in $\theta$ if $P_\theta\{\mathbf{X} \in \mathbf{C}\}$ is increasing in $\theta$ for any set $\mathbf{C}$ that is increasing.

EXAMPLE 1 (Random variables with mixtures of independent stochastically increasing distributions). In multiple testing, the $X_i$'s often have distributions that are mixtures of independent stochastically increasing distributions. That is, the density of $P_\theta$ is often of the form

$$f_\theta(\mathbf{x}) = \int \prod_{i=1}^n f_{i\theta_i}(x_i, y)\, dG(y),$$

where $f_{i\theta_i}(x, y)$ is stochastically increasing in $\theta_i$ for each $y$ and $G$ is a probability distribution independent of $\theta$. A stronger condition—which is that for any $\theta_i < \theta_i'$, $f_{i\theta_i'}(x, y)/f_{i\theta_i}(x, y)$ is increasing in $x$ for each $i$, the monotone likelihood ratio (MLR) condition of Lehmann [12] satisfied by many of the commonly used distributions—is often useful to check for the stochastic increasing property of $f_{i\theta_i}(x, y)$ in $\theta_i$. The multivariate distribution of such random variables is stochastically increasing in $\theta$.

EXAMPLE 2 (Multivariate location family of distributions). Let the density of $P_\theta$ be of the form $f_\theta(\mathbf{x}) \equiv f(\mathbf{x} - \theta)$. Distributions of this type are stochastically increasing. This is because, for any $\theta < \theta'$, we have

$$P_{\theta'}\{\mathbf{X} \in C\} = P_\theta\{\mathbf{X} \in C - (\theta' - \theta)\} \geq P_\theta\{\mathbf{X} \in C\}.$$

Many of the distributions that arise in multiple testing are of the type in Example 1 or 2. For instance, (i) independent normals with $\theta_i$'s representing the means, (ii) absolute values of independent normals with $\theta_i$'s representing the absolute means, (iii) independent chi-squares where $\theta_i$'s are the scale parameters or (iv) scaled mixtures of all these distributions, are of the type in Example 1. They arise in simultaneous testing of means or variances of independent normals against one- or two-sided alternatives. Multivariate $\ln F$ that arises in many-to-one comparisons of variances against one-sided alternatives is another distribution of the type in Example 1. Multivariate normal and multivariate $t$ are distributions of the type in Example 2, arising, for instance, in Dunnett's many-to-one comparisons of means against one-sided alternatives in a one-way layout with a known or unknown common variance.

**3. Results on FDR for fixed $J_0$ and $J_1$.** In this section we derive results on the FDR of a single-step procedure, assuming fixed, but unknown, $J_0$ and $J_1$. We use the following notation here and in the rest of the paper. Define $J = \{1, \ldots, n\}$ and $J_{(-i)} = J - \{i\}$. Define $X_{(1)} \leq \cdots \leq X_{(n)}$ as the ordered components of the set $\{X_j : j \in J\}$ and $X_{(1)}^{(-i)} \leq \cdots \leq X_{(n-1)}^{(-i)}$ as those of the subset $\{X_j : j \in J_{(-i)}\}$. We assume that the marginal distribution of any $X_i$ depends on $\theta$ only through the corresponding $\theta_i$.

First, we have the following lemma.

LEMMA 1.    *The FDR of the single-step procedure with fixed critical value $t$ is given by*

$$
\begin{aligned}
\text{FDR}_\theta & (t; J_0, J_1) \\
(3.1) \quad &= \sum_{i \in J_0} \left[ P_{\theta_i}\{X_i \geq t\} - \sum_{j=1}^{n-1} \frac{P_\theta\{X_{(j)}^{(-i)} \geq t, X_i \geq t\}}{(n-j)(n-j+1)} \right] \\
&= P_\theta\{X_{(n)} \geq t\} - \sum_{i \in J_1} \left[ P_{\theta_i}\{X_i \geq t\} - \sum_{j=1}^{n-1} \frac{P_\theta\{X_{(j)}^{(-i)} \geq t, X_i \geq t\}}{(n-j)(n-j+1)} \right].
\end{aligned}
$$

Now suppose that $\mathbf{X}$ is stochastically increasing in $\theta$. Then, since the set $\{X_{(j)}^{(-i)} \geq t, X_i \geq t\}$ is increasing in $\mathbf{X}$, the probability $P_\theta\{X_{(j)}^{(-i)} \geq t, X_i \geq t\}$ is increasing in $\theta$. The probability $P_\theta\{X_{(n)} \geq t\}$ is also increasing in $\theta$ because $\{X_{(n)} \geq t\}$ is an increasing set. Thus, using the first expression of the FDR in (3.1), we notice that it is decreasing in $\theta$ and, hence, in $\{\theta_i : i \in J_1\}$ for fixed $\{\theta_i : i \in J_0\}$, whereas from the the second expression we see that it is increasing in $\{\theta_i : i \in J_0\}$ for fixed $\{\theta_i : i \in J_1\}$. In other words, $\text{FDR}_\theta(t; J_0, J_1)$ decreases as $\theta_i$ moves away from $\theta_{i0}$ for at least one $i \in J_0$ or at least one $i \in J_1$, with

$$
(3.2) \qquad\qquad \sup_\theta \text{FDR}_\theta(t; J_0, J_1) = \text{FDR}_{\theta_0}(t; J_0, J_1),
$$

where $\theta_0 = (\theta_{10}, \ldots, \theta_{n0})$. If $\mathbf{X}$ is exchangeable when $\theta = \theta_0$ with the common marginal c.d.f. $F_0$, the right-hand side of (3.2) reduces to

$$
\begin{aligned}
(3.3) \quad & n_0 \left[ \bar{F}_0(t) - \sum_{j=1}^{n-1} \frac{P_{\theta_0}\{X_{(j)}^{(-1)} \geq t, X_1 \geq t\}}{(n-j)(n-j+1)} \right] \\
&= \frac{n_0}{n} \text{FDR}_{\theta_0}(t; J, \phi) \\
&= \frac{n_0}{n} P_{\theta_0}\{R > 0\},
\end{aligned}
$$

where $\bar{F}_0 = 1 - F_0$ and $\phi$ represents a null set. Thus, we have the following theorem, which is one of the main results of this article.

THEOREM 1. *If* **X** *is stochastically increasing in* $\theta$, *then* $\mathrm{FDR}_\theta(t, J_0, J_1)$ *decreases as* $\theta_i$ *moves away from* $\theta_{i0}$ *for at least one* $i \in J_0$ *or for at least one* $i \in J_1$. *Furthermore, if* **X** *is exchangeable when* $\theta = \theta_0$, *then*

$$(3.4) \qquad \sup_\theta \mathrm{FDR}_\theta(t; J_0, J_1) = \frac{n_0}{n} P_{\theta_0}\{R > 0\}.$$

Theorem 5.3 of [5] gives the above decreasing property of FDR with respect to only $\{\theta_i, i \in J_1\}$ under the assumptions that $\{X_i, i \in J_0\}$ and $\{X_i, i \in J_1\}$ are jointly independent and $\{X_i, i \in J_1\}$ is stochastically increasing in $\{\theta_i, i \in J_1\}$. Theorem 1 is a version of this for single-step procedures with dependent **X** and one-sided null hypotheses.

As a corollary to Theorem 1, if the critical value $t$ provides a level $\alpha$ test for the overall null hypothesis $\bigcap_{i=1}^n H_i$, that is, if $t$ satisfies $P_{\theta_0}\{R > 0\} = P_{\theta_0}\{\max_{i \in J} X_i \geq t\} \leq \alpha$, then we have

$$(3.5) \qquad \mathrm{FDR}_\theta(t; J_0, J_1) \leq \frac{n_0}{n}\alpha,$$

implying that the FDR is controlled at $\alpha$. Inequality (3.5) is interesting in that it represents a single-step analog of the same inequality known to hold for stepwise procedures with Simes [16] critical values providing an $\alpha$-level test for $\bigcap_{i=1}^n H_i$ [1, 5, 14]. Regarding the choice for $t$, if one does not want to utilize the distributional form of **X** or if it is unknown, the Bonferroni critical value that satisfies

$$(3.6) \qquad F_0(t) = 1 - \frac{\alpha}{n}$$

can be used. If, however, **X** is known to be positively dependent so that the inequality $P_{\theta_0}\{\max_{i \in J} X_i < t\} \geq F_0^n(t)$ holds under the null hypothesis values with the equality holding under independence, as in the case of many distributions that arise in multiple testing, the Šidák critical value $t$ that satisfies the equation

$$(3.7) \qquad F_0(t) = (1 - \alpha)^{1/n}$$

offers a less conservative choice.

We should point out that there is no surprise that the Bonferroni and Šidák single-step procedures control FDR, because they are known to control the familywise error rate (FWER). It is also known that, given $n_0$, it can be incorporated in the Bonferroni and other procedures to improve their FWER control [10]. What is new here is that Bonferroni and Šidák procedures can be further improved in terms of having better control of FDR using an estimate of $n_0$, in the spirit of Benjamini and Hochberg [2], Benjamini, Krieger and Yekutieli [3], Storey [17] and Storey, Taylor and Siegmund [19]. For instance, since $\sup_\theta \mathrm{FDR}_\theta(t; J_0, J_1) \leq n_0\{1 - F_0(t)\}$, as we see from Theorem 1, rather than controlling $n\{1 - F_0(t)\}$, which the Bonferroni method does, a better control of FDR can be achieved if we control $\hat{n}_0\{1 - F_0(t)\}$ for some appropriately chosen estimate $\hat{n}_0$ of $n_0$. To estimate $n_0$, Storey [17] suggested using the ratio $K_\tau/F_0(\tau)$, where $K_\tau = \sum_{i=1}^n I(X_i < \tau)$, for some well-chosen $\tau$. However, Storey, Taylor and Siegmund [19] slightly modified it and used

$$(3.8) \qquad \hat{n}_0(\tau) = \frac{K_\tau + 1}{F_0(\tau)}$$

to obtain a new class of BH-type FDR-controlling procedures under independence. We use this $\hat{n}_0$ in our modification to the Bonferroni procedure. Also, the $X_i$'s that are small compared to $\tau$ should not be declared large when modified Bonferroni is used. Thus, our modified Bonferroni procedure rejects $H_i$ whenever

$$(3.9) \qquad X_i \geq \max\left\{\tau, F_0^{-1}\left(1 - \frac{\alpha F_0(\tau)}{K_\tau + 1}\right)\right\}.$$

We prove later in this section that our modified Bonferroni procedure controls FDR under independence and we provide numerical evidence showing that quite often this control can be achieved much less conservatively. However, when $\mathbf{X}$ is known to be independent or at least positively dependent, a modification to the Šidák procedure is expected to produce a better performing procedure than the modified Bonferroni procedure. So, we first modify the Šidák procedure. The following theorem suggests how the idea of modifying the Bonferroni procedure can be extended to that for the Šidák procedure. It extends the inequality for the FDR under independence, given by Theorem 1, from a single-step to a two-step procedure that, for some fixed $\tau \in (-\infty, \infty)$ and a predetermined function $t_\tau(k) \geq \tau$, $k = 0, 1, \ldots, n$, first finds $k = \max_{0 \leq i \leq n}\{i : X_{(i)} < \tau\}$ (note that $X_{(0)} = -\infty$), then rejects all $H_i$ for which $X_i \geq t_\tau(k)$.

THEOREM 2. *Let $\mathbf{X}$ be independent with the distribution of $X_i$, indexed by the parameter $\theta_i$, belonging to an MLR family and having identical marginals when $\theta = \theta_0$. Then, for a two-step procedure with $t_\tau(k) \geq \tau$, for all $k = 0, 1, \ldots, n$, the FDR satisfies the inequality*

$$\begin{aligned}(3.10) \quad &\mathrm{FDR}_\theta^{(2)}(t_\tau \geq \tau; J_0, J_1) \\ &\leq \bar{F}_0(\tau) \sum_{i \in J_0} \sum_{k=0}^{n-1} \frac{1}{n-k}\left[1 - \left(1 - \frac{\bar{F}_0(t_\tau(k))}{\bar{F}_0(\tau)}\right)^{n-k}\right] \\ &\qquad \times P_\theta\{X_{(k)}^{(-i)} < \tau \leq X_{(k+1)}^{(-i)}\}\end{aligned}$$

*(with $X_{(0)}^{(-i)} = -\infty$ and $X_{(n)}^{(-i)} = \infty$).*

When $\tau = -\infty$, $k = 0$ with probability 1 and (3.10) reduces to the one given by Theorem 1 under independence with $t = t_{-\infty}(0)$. It is interesting to see that $\text{FDR}_\theta^{(2)}(t_\tau \geq \tau; J_0, J_1) \leq \text{FDR}_\theta(\tau; J_0, J_1)$.

The modified Bonferroni procedure is a two-step procedure with $t_\tau(k)$ given by the right-hand side of (3.9) given $K_\tau = k$; that is, $t_\tau(k)$ is such that $\bar{F}_0(t_\tau(k)) = \min\{\bar{F}_0(\tau), \alpha F_0(\tau)/(k+1)\}$. We propose to modify the Šidák procedure using a two-step procedure where $t_\tau(k)$ is such that

$$(3.11) \quad \bar{F}_0(t_\tau(k)) = \bar{F}_0(\tau)\left[1 - \left(1 - \min\left\{1, \frac{\alpha(n-k)F_0(\tau)}{(k+1)\bar{F}_0(\tau)}\right\}\right)^{1/(n-k)}\right]$$

with $t_\tau(n) = \infty$. The right-hand side of (3.10) for this modified Šidák procedure is less than or equal to

$$(3.12) \quad \begin{aligned} &\alpha \sum_{i \in J_0} \sum_{k=0}^{n-1} \frac{F_0(\tau)}{k+1} P_\theta\{X_{(k)}^{(-i)} < \tau \leq X_{(k+1)}^{(-i)}\} \\ &\leq \alpha \sum_{i \in J} \sum_{k=1}^{n} \frac{1}{k} P_\theta\{X_i < \tau, X_{(k-1)}^{(-i)} < \tau \leq X_{(k)}^{(-i)}\} \\ &= \alpha \sum_{k=1}^{n} P_\theta\{X_{(k)} < \tau \leq X_{(k+1)}\} \\ &= \alpha P_\theta\{X_{(1)} < \tau\}; \end{aligned}$$

see, for example, [13], page 497, for the first equality in (3.12). Thus, we see that our modified Šidák procedure controls FDR under independence.

The right-hand side of (3.10) is less than or equal to

$$(3.13) \quad \sum_{i \in J_0} \sum_{k=0}^{n-1} \bar{F}_0(t_\tau(k)) P_\theta\{X_{(k)}^{(-i)} < \tau \leq X_{(k+1)}^{(-i)}\},$$

which, for the modified Bonferroni procedure, is less than or equal to the first expression in (3.12). Thus, the FDR of the modified Bonferroni procedure is also less than or equal to $\alpha P_\theta\{X_{(1)} < \tau\}$ and, hence, is controlled; of course, it is controlled more conservatively than the modified Šidák procedure.

We conducted a numerical study to investigate the extent of improvement offered by our modified Šidák procedure in controlling FDR over the modified Bonferroni and the original Bonferroni and Šidák procedures. We generated $n = 100$ dependent random variables $X_i \sim N(\mu_i, 1)$, $i = 1, \ldots, 100$, with the same variance 1 and a common correlation $\rho$, and performed 100 hypothesis tests of $\mu = 0$ against $\mu > 0$, each using first the Bonferroni critical value and then the Šidák critical value corresponding to $\alpha = 0.05$. The value of $Q$ was then calculated for each procedure by setting $n_0$ of the $\mu_i$'s to zero and the remaining $\mu_i$'s to a positive value $\delta$. The FDR then was estimated by averaging the $Q$ values over 5000 iterations. Thus, we have the simulated

TABLE 2
*Simulated values of the FDR of the Bonferroni and Šidák procedures and their modifications with $\alpha = 0.05$*

| | | Independent ($\rho = 0$) | | | | Dependent ($\rho = 0.5$) | | | |
|---|---|---|---|---|---|---|---|---|---|
| | | Bonferroni | | Šidák | | Bonferroni | | Šidák | |
| $n_0$ | $\delta$ | Original | Modified | Original | Modified | Original | Modified | Original | Modified |
| 30 | 0.5 | 0.0118 | 0.0150 | 0.0119 | 0.0167 | 0.0048 | 0.0167 | 0.0049 | 0.0332 |
| | 1.5 | 0.0045 | 0.0073 | 0.0045 | 0.0079 | 0.0006 | 0.0066 | 0.0006 | 0.0412 |
| | 2.5 | 0.0008 | 0.0019 | 0.0008 | 0.0022 | 0.0002 | 0.0054 | 0.0002 | 0.0412 |
| 50 | 0.5 | 0.0218 | 0.0259 | 0.0222 | 0.0276 | 0.0092 | 0.0307 | 0.0093 | 0.0493 |
| | 1.5 | 0.0103 | 0.0147 | 0.0106 | 0.0149 | 0.0015 | 0.0116 | 0.0015 | 0.0455 |
| | 2.5 | 0.0021 | 0.0031 | 0.0021 | 0.0033 | 0.0006 | 0.0093 | 0.0006 | 0.0441 |
| 70 | 0.5 | 0.0315 | 0.0349 | 0.0319 | 0.0359 | 0.0141 | 0.0488 | 0.0144 | 0.0667 |
| | 1.5 | 0.0187 | 0.0237 | 0.0189 | 0.0232 | 0.0034 | 0.0196 | 0.0034 | 0.0494 |
| | 2.5 | 0.0052 | 0.0061 | 0.0052 | 0.0061 | 0.0013 | 0.0154 | 0.0014 | 0.0463 |
| 90 | 0.5 | 0.0414 | 0.0423 | 0.0423 | 0.0434 | 0.0234 | 0.0734 | 0.0240 | 0.0903 |
| | 1.5 | 0.0351 | 0.0382 | 0.0359 | 0.0393 | 0.0108 | 0.0414 | 0.0111 | 0.0642 |
| | 2.5 | 0.0173 | 0.0180 | 0.0175 | 0.0189 | 0.0045 | 0.0311 | 0.0046 | 0.0554 |
| | MaxSE | 0.0028 | 0.0028 | 0.0028 | 0.0029 | 0.0020 | 0.0034 | 0.0021 | 0.0038 |

FDR of the Bonferroni and Šidák procedures. We chose $F_0(\tau) = 1/2$ and similarly calculated the FDR of the modified Bonferroni and Šidák procedures corresponding to this $\tau$. Table 2 compares the FDRs of the Bonferroni and Šidák procedures and their modification for $n_0 = 30, 50, 70$ and $90$, $\rho = 0$ (independent) and 0.5 (dependent), and for different values of $\delta$. The last row of this table gives the maximum of the standard errors of the estimated (simulated) FDRs in each column.

As we expected, the modified Šidák procedure provided the least conservative control of FDR under independence. Since the Bonferroni and Šidák procedures are relatively more conservative when the actual proportion of true null hypotheses is small, the idea of improving them using an estimate of $n_0$ should work well in this situation. This idea is confirmed by our numerical study. Both modified Bonferroni and modified Šidák procedures are seen to control FDR much less conservatively than their unmodified versions under independence. In the dependent case, however, the idea of improving the Bonferroni and Šidák procedures may not work unless $n_0$ is small and the dependence is weak.

Having found more than one procedure that can control the FDR under independence (e.g., the Bonferroni, Šidák and their modifications), comparing them further in terms of *power* seems to be the next important objective. While the idea of power can be conceptualized in terms of Type II errors (false negatives) in several different ways, extending it from single testing to multiple testing, one particular concept, which is the average power [i.e., $\frac{1}{n_1}E(S)$], has been used in a number of recent papers to compare FDR-controlling procedures [4, 17, 19]. However, it is argued in [15] that since the FDR is a measure of false positives, it seems more appropriate to compare different FDR-controlling procedures using a similar measure in terms of false negatives, the FNR [8, 15]. It will be interesting to see how the different FDR-controlling procedures in this paper compare in terms of measures involving FNR under the same distributional setting. This will be carried out in Section 5 after deriving some results on FNR in the next section.

**4. Results on FNR for fixed $J_0$ and $J_1$.** We will derive in this section some results on FNR of a single-step procedure, analogous to those on FDR, again assuming a fixed configuration of true and false null hypotheses. First, we have the following lemma.

LEMMA 2. *An explicit expression of FNR is*

$$
\begin{aligned}
&\mathrm{FNR}_\theta(t; J_0, J_1) \\
&= \sum_{i \in J_1}\left[ P_{\theta_i}\{X_i < t\} - \sum_{j=1}^{n-1} \frac{P_\theta\{X_{(j)}^{(-i)} < t, X_i < t\}}{j(j+1)} \right] \\
&= P_\theta\{X_{(1)} < t\} - \sum_{i \in J_0}\left[ P_{\theta_i}\{X_i < t\} - \sum_{j=1}^{n-1} \frac{P_\theta\{X_{(j)}^{(-i)} < t, X_i < t\}}{j(j+1)} \right].
\end{aligned}
$$

(4.1)

Making the same kind of arguments as we made before for the monotonicity property of the FDR, we notice that if $\mathbf{X}$ is stochastically increasing in $\theta$, the FNR is increasing in $\{\theta_i : i \in J_0\}$ for fixed $\{\theta_i : i \in J_1\}$ and is decreasing in $\{\theta_i : i \in J_1\}$ for fixed $\{\theta_i : i \in J_0\}$. In other words, $\mathrm{FNR}_\theta(t; J_0, J_1)$ decreases as $\theta_i$ moves away from $\theta_{i0}$ for at least one $i \in J_0$ or at least one $i \in J_1$, with

$$
\sup_\theta \mathrm{FNR}_\theta(t; J_0, J_1) = \mathrm{FNR}_{\theta_0}(t; J_0, J_1). \tag{4.2}
$$

Since, when $\theta = \theta_0$, $\mathbf{X}$ is exchangeable, the right-hand side in (4.2) reduces to

$$
n_1\left[ F_0(t) - \sum_{j=1}^{n-1} \frac{P_{\theta_0}\{X_{(j)}^{(-1)} < t, X_1 < t\}}{j(j+1)} \right] = \frac{n_1}{n} P_{\theta_0}\{A > 0\}. \tag{4.3}
$$

The equality in (4.3) follows from (4.1); see also [13]. This gives the next main result of this article.

THEOREM 3. *If* **X** *is stochastically increasing in* $\theta$, *then* $\mathrm{FNR}_\theta(t, J_0, J_1)$ *decreases as* $\theta_i$ *moves away from* $\theta_{i0}$ *for at least one* $i \in J_0$ *or for at least one* $i \in J_1$. *Furthermore, if* **X** *is exchangeable when* $\theta = \theta_0$, *then*

$$(4.4) \qquad \sup_\theta \mathrm{FNR}_\theta(t; J_0, J_1) = \frac{n_1}{n} P_{\theta_0}\{A > 0\}.$$

Clearly, the FNR of a single-step procedure can be controlled at a level $\beta$ under the condition stated in the above theorem by choosing a fixed $t$ subject to the condition $P_{\theta_0}\{A > 0\} = P_{\theta_0}\{\min_{i \in J} X_i \leq t\} \leq \beta$. If the dependence structure of **X** is not utilized, the equation $F_0(t) = \beta/n$ provides a Bonferroni-type choice for $t$. When **X** is known to be positively dependent so that the inequality $P_{\theta_0}\{\min_{i \in J} X_i \geq t\} \geq \bar{F}_0^n(t)$ is true, with the equality holding under independence, Šidák-type $t$ can be determined from the equation $F_0(t) = 1 - (1 - \beta)^{1/n}$. These procedures can potentially be improved in terms of having better control of FNR by borrowing information from the $X_i$'s exceeding an appropriately chosen value $\tau$.

The following theorem is a FNR analog of Theorem 2 that extends the inequality on FNR given by Theorem 3 from a single-step to a two-step procedure and suggests how to modify the above single-step FNR-controlling procedures.

THEOREM 4. *Under the conditions stated in Theorem* 2, *the FNR of a two-step procedure with* $t_\tau(k) \leq \tau$ *for all* $k = 0, 1, \ldots, n$ *satisfies the inequality*

$$
\begin{aligned}
\mathrm{FNR}_\theta^{(2)}&(t_\tau \leq \tau; J_0, J_1) \\
&\leq F_0(\tau) \sum_{i \in J_1} \sum_{k=1}^n \frac{1}{k} \left[ 1 - \left( 1 - \frac{F_0(t_\tau(k))}{F_0(\tau)} \right)^k \right] P_\theta\{X_{(k-1)}^{(-i)} < \tau \leq X_{(k)}^{(-i)}\}.
\end{aligned}
$$
$(4.5)$

When $\tau \to \infty$, $k = n$ with probability 1 and the above inequality reduces to that given by Theorem 3 under independence with $t = t_\infty(n)$. We modify the Šidák procedure using a two-step procedure with $t_\tau(k) \leq \tau$ satisfying

$$(4.6) \quad F_0(t_\tau(k)) = F_0(\tau) \left[ 1 - \left( 1 - \min\left\{ 1, \frac{\beta k \bar{F}_0(\tau)}{(n-k+1)F_0(\tau)} \right\} \right)^{1/k} \right]$$

and $t_\tau(0) = -\infty$. For this modified Šidák procedure,

$$
\begin{aligned}
\mathrm{FNR}_\theta^{(2)}&(t_\tau \leq \tau; J_0, J_1) \\
&\leq \beta \sum_{i \in J_1} \sum_{k=1}^n \frac{\bar{F}_0(\tau)}{n-k+1} P_\theta\{X_{(k-1)}^{(-i)} < \tau \leq X_{(k)}^{(-i)}\} \\
&\leq \beta \sum_{i \in J_1} \sum_{k=0}^{n-1} \frac{1}{n-k} P_\theta\{X_i \geq \tau, X_{(k)}^{(-i)} < \tau \leq X_{(k+1)}^{(-i)}\} \\
&\leq \beta \sum_{i \in J} \sum_{k=0}^{n-1} \frac{1}{n-k} P_\theta\{X_i \geq \tau, X_{(k)}^{(-i)} < \tau \leq X_{(k+1)}^{(-i)}\} \\
&= \beta \sum_{k=0}^{n-1} P_\theta\{X_{(k)} < \tau \leq X_{(k+1)}\} \\
&= \beta P_\theta\{X_{(n)} \geq \tau\}.
\end{aligned}
\tag{4.7}
$$

The second inequality in (4.7) follows from the fact that $\bar{F}_0(\tau) \leq \bar{F}_{\theta_i}(\tau)$; for the first equality, see [13]. Thus, the above modified Šidák procedure controls FNR under independence.

The right-hand side of (4.5) is less than or equal to

$$
\sum_{i \in J_1} \sum_{k=1}^n F_0(t_\tau(k)) P_\theta\{X_{(k-1)}^{(-i)} < \tau \leq X_{(k)}^{(-i)}\}.
\tag{4.8}
$$

This is less than or equal to the right-hand side of the first inequality in (4.7), which is less than or equal to $\beta$, if we choose $t_\tau(k) \leq \tau$ satisfying

$$
F_0(t_\tau(k)) = \min\left\{F_0(\tau), \frac{\beta \bar{F}_0(\tau)}{n-k+1}\right\}.
\tag{4.9}
$$

This gives us our FNR-controlling modified Bonferroni procedure, which is of course more conservative than the modified Šidák procedure in the sense that it allows less nondiscoveries.

REMARK 1. It is important to note that the above results on FNR have been developed with the idea of controlling false nondiscoveries of any set of true alternatives (or false nulls). However, one is often interested in controlling false nondiscoveries of a prespecified set of true alternatives. These results can be easily modified in such a situation. Let $\theta_i = \theta_{i1}$ for some specified $\theta_{i1} > \theta_{i0}$, $i \in J_1$. Assume that $\mathbf{X}$ is exchangeable under $\theta = \theta_1 = (\theta_{11}, \ldots, \theta_{n1})$. Then Theorem 3 can be modified to

$$
\sup_\theta \mathrm{FNR}(t; J_0, J_1) = \frac{n_1}{n} P_{\theta_1}\{A > 0\}
\tag{4.10}
$$

and Theorem 4 can be modified to

$$
\begin{aligned}
&\mathrm{FNR}_\theta^{(2)}(t_\tau \leq \tau; J_0, J_1)\\
(4.11) \quad &\leq F_1(\tau) \sum_{i \in J_1} \sum_{k=1}^{n} \frac{1}{k}\left[1 - \left(1 - \frac{F_1(t_\tau(k))}{F_1(\tau)}\right)^k\right]\\
&\qquad\qquad \times P_\theta\{X_{(k-1)}^{(-i)} < \tau \leq X_{(k)}^{(-i)}\},
\end{aligned}
$$

where $F_1$ is the common c.d.f. of $X_i$ under $\theta_{i1}$. The Bonferroni and Šidák procedures as well as their two-step modifications using critical values based on $F_1$ will provide better control of FNR in this case than values based on $F_0$.

We conducted a numerical study to investigate how well these different FNR procedures control FNR under a specified set of true alternatives. We noticed, as in the case of controlling FDR, that although both modified Bonferroni and Šidák procedures often control FNR much less conservatively than their unmodified versions, the modified Šidák procedure provides the best control of FNR.

**5. A numerical study.** In this section we compare the different FDR-controlling procedures under independence discussed in Section 3 in terms of a concept of power that relates to the unbiasedness condition Sarkar [15] introduced. Since the FDR measures the expected proportion of incorrect decisions, a good multiple testing procedure must ensure that it does not exceed the expected proportion of correct decisions. The quantity $1 - \mathrm{FNR}$, which Genovese and Wasserman [8] called the correct nondiscovery rate, is a measure of correct decisions. In situations where controlling false negatives is of primary importance, the FNR provides a measure of incorrect decisions with the corresponding measure of correct decisions being $1 - \mathrm{FDR}$. Whether we have a multiple testing procedure designed to control FDR or FNR, the inequality $\mathrm{FDR} + \mathrm{FNR} \leq 1$ represents a desirable property for any such multiple testing procedure. This is referred to as the unbiasedness condition of an FDR- or FNR-controlling multiple testing procedure. A natural way to compare different FDR- or FNR-controlling procedures would be to see how they perform in terms of a measure that reflects the strength of unbiasedness. This leads us to the consideration of the quantity

$$
(5.1) \qquad\qquad \pi_\theta = 1 - \mathrm{FDR}_\theta - \mathrm{FNR}_\theta.
$$

It is also related to the idea of Genovese and Wasserman [8], who suggested using $1 - \pi_\theta$ as a risk function to compare multiple testing procedures. This is our concept of power.

We investigated how the different FDR procedures in Section 3 perform in terms of the aforementioned concept of power. We computed the FNR

and then the power $1 - \text{FNR} - \text{FDR}$ for the Bonferroni and Šidák procedures and their modified versions [with $F_0(\tau) = 1/2$] based on the normal data that have been simulated before for FDR calculations. These simulated powers are displayed in Figure 1. As we see from this figure, the modified Šidák procedure is often the most powerful under independence, especially, as one would expect, when the proportion of true null hypotheses is relatively small. The unmodified Bonferroni and Šidák procedures, not surprisingly, are practically indistinguishable in terms of their power performance. One should, however, be cautious in interpreting this graph in the dependent case (particularly, the upper right two panels), in light of Table 1, which indicates that the modified Bonferroni and Šidák procedures may fail to control FDR unless the dependence is weak and $n_0$ is small.

We should point out that the unbiasedness property of the single-step procedures, which is numerically seen to hold, can be theoretically proved easily from Theorems 1 and 3. However, a theoretical justification of the same property for the two-step procedures, which appears to be also true from Figure 1, is an interesting and a more challenging theoretical problem.

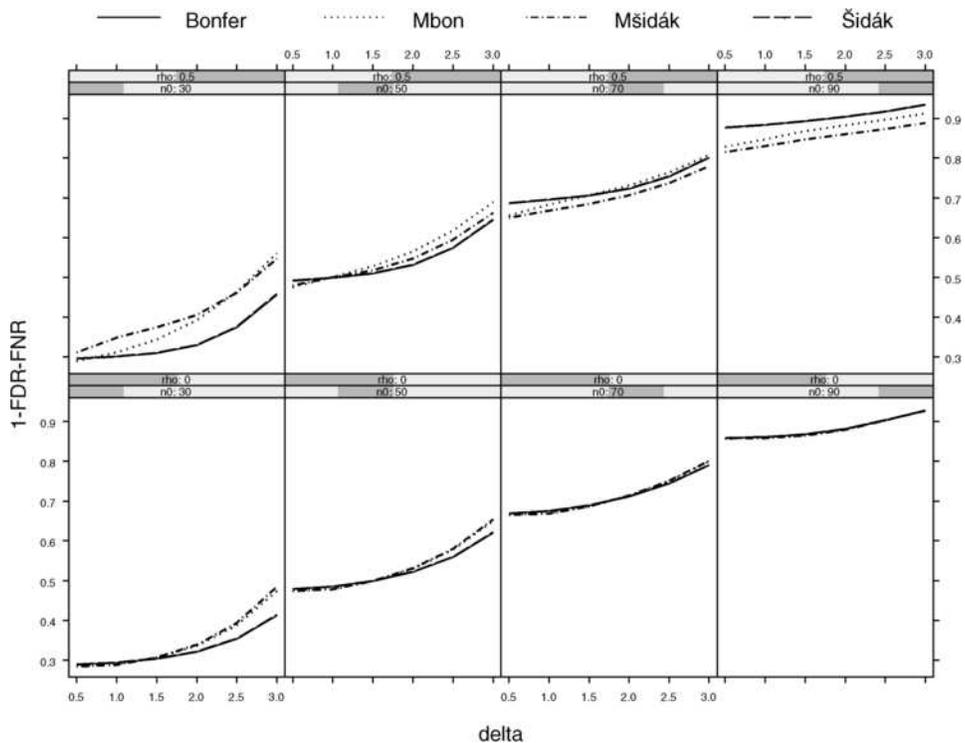

FIG. 1. *Comparison of Bonferroni and Šidák procedures with their modified versions in terms of* $1 - FDR - FNR$.

Also, the same concept of power could be used to compare different FNR-controlling procedures.

**6. Results on FDR and FNR under a mixture model.** In this section, we present appropriate modifications to Lemmas 1 and 2 when a mixture approach is taken as in [7, 17]. We will, however, assume a slightly more general mixture model in the sense that it does not assume independence of the test statistics. More specifically, we first let $\mathbf{H} = (H_1, \ldots, H_n)$, with $H_i = 0$ indicating that $H_i$ is true and $H_i = 1$ indicating that it is false. Then we assume that $(X_i, H_i)$, $i = 1, \ldots, n$, have the distribution

$$\mathbf{X}|\mathbf{H} \sim f(\mathbf{x}, \theta_{\mathbf{H}}) \qquad \text{where } \theta_{\mathbf{H}} = (\theta_{H_1}, \ldots, \theta_{H_n}), \theta_{H_i} = (1 - H_i)\theta_i' + H_i\theta_i'',$$
$$\text{with } \theta_i' \leq \theta_{i0}, \theta_i'' > \theta_{i0}, i = 1, \ldots, n,$$

and $\mathbf{H} \sim \pi_{\mathbf{h}}$, where $\pi_{\mathbf{h}}$ are some probabilities defined on
$\mathcal{H} = \{\mathbf{h} = (h_1, \ldots, h_n) : h_i = 0 \text{ or } 1\}$.

(6.1)

Regarding $f$, we assume that it belongs to a location family of distributions; that is, $f(\mathbf{x}, \theta_{\mathbf{H}}) = f(\mathbf{x} - \theta_{\mathbf{H}})$, with a positive dependence structure that ensures that, for any increasing (or decreasing) function $\phi$ of $\mathbf{X}$, the expectation $E\{\phi(\mathbf{X})|X_i, \mathbf{H}\}$ is increasing (or decreasing) in $X_i$. This is true if, for instance, $\mathbf{X}$ is positive regression dependent on subset (PRDS) under the density $f(\mathbf{x})$, as in the case of multivariate normal with positive correlations and many other multivariate distributions encountered in multiple testing; see, for example, [5, 14]. Of course, when $(X_i, H_i)$, $i = 1, \ldots, n$, are independent, we assume no particular form for the density $f$; that is, we simply assume that $X_i|H_i \sim f(x, \theta_{H_i})$. Since we assume that $\theta_i$ takes the value $\theta_i'$ when $H_i = 0$ and the value $\theta_i''$ when $H_i = 1$, the probabilities in the following discussion are all evaluated under these fixed $\theta' = (\theta_1', \ldots, \theta_n')$ and $\theta'' = (\theta_1'', \ldots, \theta_n'')$.

THEOREM 5. *Under the above mixture model and the conditions assumed therein,*

$$\text{FDR}(t, n) \leq \sum_{i=1}^{n} \delta_i P\{H_i = 0 | X_i \geq t\},$$

(6.2)

*where*

$$\delta_i = P\{X_i \geq t\} - \sum_{j=1}^{n-1} \frac{P\{X_{(j)}^{(-i)} \geq t, X_i \geq t\}}{(n-j)(n-j+1)} \quad \text{and}$$
$$\sum_{i=1}^{n} \delta_i = P\{R > 0\},$$

(6.3)

*with the equality holding when the $(X_i, H_i)$'s are independent.*

When $(X_i, H_i)$, $i = 1, \ldots, n$, are identically distributed, Theorem 5 reduces to

$$\text{FDR}(t, n) \leq P\{H_1 = 0 | X_1 \geq t\} P\{R > 0\}. \tag{6.4}$$

The equality in (6.4) holds when $(X_i, H_i)$, $i = 1, \ldots, n$, are i.i.d., which is Storey's [17, 18] result, providing a "Bayesian Type I error rate" interpretation to his notion of pFDR = FDR/$P\{R > 0\}$. Thus, the following corollary to Theorem 5 is an extension of his result to the dependent case.

COROLLARY 1.   *Under the above mixture model and the conditions assumed therein,*

$$\text{pFDR}(t, n) \leq \max_{1 \leq i \leq n} P\{H_i = 0 | X_i \geq t\}. \tag{6.5}$$

*When the $(X_i, H_i)$'s are identically distributed, we have*

$$\text{pFDR}(t, n) \leq P\{H_1 = 0 | X_1 \geq t\}, \tag{6.6}$$

*with the equality holding when the $(X_i, H_i)$'s are i.i.d.*

Storey [17] introduced a pFDR analog of the $p$-value, called the $q$-value, that provides a measure of the strength of the tests in a multiple testing procedure with respect to pFDR. For a single-step multiple testing procedure of $n$ hypotheses with a rejection region of the form $X_i \geq t$ for each $H_i$, it is defined as

$$q_n(t) = \inf_{x \leq t} \text{pFDR}(x, n). \tag{6.7}$$

Storey [17], however, considered this quantity when $(X_i, H_i)$, $i = 1, \ldots, n$, are i.i.d., which is

$$q(t, H_1) = \inf_{x \leq t} P\{H_1 = 0 | X_1 \geq x\}. \tag{6.8}$$

Corollary 1 says that when the $(X_i, H_i)$'s are dependent with common marginals, in the sense assumed in that corollary, we have $q_n(t) \leq q(t, H_1)$. That is, the $q$-value of a single-step multiple test procedure obtained under certain commonly encountered types of dependence is more conservative, as one would want, compared to the corresponding i.i.d. case.

THEOREM 6.   *Under the conditions stated in Theorem 5,*

$$\text{FNR} \leq \sum_{i=1}^{n} \gamma_i P\{H_i = 1 | X_i < t\}, \tag{6.9}$$

*where*

$$(6.10) \qquad \gamma_i = P\{X_i < t\} - \sum_{j=1}^{n-1} \frac{P\{X_{(j)}^{(-i)} < t, X_i < t\}}{j(j+1)} \quad and$$

$$\sum_{i=1}^{n} \gamma_i = P\{A > 0\},$$

*with the equality holding when the* $(X_i, H_i)$*'s are independent.*

This theorem can be proved following arguments similar to those used to prove Theorem 5 and with the help of an identity for $P(A > 0)$ given by Sarkar [13].

COROLLARY 2.    *Under the conditions stated in Theorem* 5,

$$(6.11) \qquad \mathrm{pFNR} \leq \max_{1 \leq i \leq n} P\{H_i = 1 | X_i < t\}.$$

*When the* $(X_i, H_i)$*'s are identically distributed, we have*

$$(6.12) \qquad \mathrm{pFNR} \leq P\{H_1 = 1 | X_1 < t\},$$

*with the equality holding when the* $(X_i, H_i)$*'s are i.i.d.*

## 7. Proofs.

PROOF OF LEMMA 1.    The FDP is given by

$$(7.1) \qquad Q(t; J_0, J_1) = \sum_{i \in J_0} \sum_{j=0}^{n-1} \frac{1}{n-j} I\{R = n-j, X_i \geq t\}.$$

Since $\{R = n - j\} = \{X_{(j)} < t \leq X_{(j+1)}\}$, with $X_{(0)} = -\infty$ and $X_{(n+1)} = \infty$, we have

$$\{R = n-j, X_i \geq t\} = \{X_{(j)}^{(-i)} < t \leq X_{(j+1)}^{(-i)}, X_i \geq t\}.$$

Therefore,

$$
\begin{aligned}
Q(t; J_0, J_1) &= \sum_{i \in J_0} \sum_{j=0}^{n-1} \frac{1}{n-j} I\{X_{(j)} < t \leq X_{(j+1)}, X_i \geq t\} \\
(7.2) \qquad &= \sum_{i \in J_0} \sum_{j=0}^{n-1} \frac{1}{n-j} [I\{X_{(j+1)}^{(-i)} \geq t, X_i \geq t\} - I\{X_{(j)}^{(-i)} \geq t, X_i \geq t\}] \\
&= \sum_{i \in J_0} I\{X_i \geq t\} - \sum_{i \in J_0} \sum_{j=1}^{n-1} \frac{I\{X_{(j)}^{(-i)} \geq t, X_i \geq t\}}{(n-j)(n-j+1)}.
\end{aligned}
$$

Taking the expectation in (7.2), we get the first expression of the FDR in Lemma 1. The second expression follows from the fact that $Q$ reduces to $I\{R > 0\} = I\{X_{(n)} \geq t\}$ if we consider the first summation in (7.2) over all $i \in J$. $\square$

PROOF OF THEOREM 2. First note that

$$
\begin{aligned}
(7.3) \quad & \mathrm{FDR}_\theta^{(2)}(t_\tau \geq \tau; J_0, J_1) \\
& = \sum_{k=0}^n E_\theta\{Q(t_\tau(k); J_0, J_1)I\{X_{(k)} < \tau \leq X_{(k+1)}\}\}.
\end{aligned}
$$

Since $t_\tau(k) \geq \tau$ for all $k$, when $k = n$ (i.e., when $X_{(n)} < \tau$), there is no rejection of null hypotheses, implying that $Q = 0$.

Let $F_{\theta_i}(x)$ and $f_{\theta_i}(x)$, respectively, be the c.d.f. and the density of $X_i$ under any alternative $\theta_i$ for $i = 1, \ldots, n$. Since the $X_i$'s are assumed to be independent, the conditional expectation of $Q(t_\tau(k); J_0, J_1)$, given $\{X_{(k)} < \tau \leq X_{(k+1)}\}$ for $k = 0, 1, \ldots, n-1$, is the FDR of the single-step procedure based on $n - k$ independent random variables $Y_1, \ldots, Y_{n-k}$ with $Y_i \sim f_{\theta_i}(x)I(x \geq \tau)/\bar{F}_{\theta_i}(\tau)$ and critical value $t_\tau(k)$. Since the density of $Y_i$ has the MLR property, implying that $(Y_1, \ldots, Y_{n-k})$ is stochastically increasing, we have from Theorem 1 that this conditional expectation is

$$
\begin{aligned}
(7.4) \quad & \leq \frac{n_0(\tau)}{n-k} P_{\theta_0}\left\{\max_{1 \leq j \leq n-k} Y_j \geq t_\tau(k)\right\} \\
& = \frac{n_0(\tau)}{n-k}\left[1 - \left(1 - \frac{\bar{F}_0(t_\tau(k))}{\bar{F}_0(\tau)}\right)^{n-k}\right],
\end{aligned}
$$

where $n_0(\tau) = \sum_{i \in J_0} I(X_i \geq \tau)$. Going back to (7.3), we then have

$$
\begin{aligned}
(7.5) \quad & \mathrm{FDR}_\theta^{(2)}(t_\tau \geq \tau; J_0, J_1) \\
& \leq \sum_{k=0}^{n-1} \sum_{i \in J_0} E_\theta\left\{\frac{1}{n-k}\left[1 - \left(1 - \frac{\bar{F}_0(t_\tau(k))}{\bar{F}_0(\tau)}\right)^{n-k}\right]\right. \\
& \qquad\qquad \left. \times I\{X_i > \tau, X_{(k)} < \tau \leq X_{(k+1)}\}\right\} \\
& = \sum_{i \in J_0} \sum_{k=0}^{n-1} \frac{1}{n-k}\left[1 - \left(1 - \frac{\bar{F}_0(t_\tau(k))}{\bar{F}_0(\tau)}\right)^{n-k}\right] \\
& \qquad\qquad \times P_\theta\{X_i > \tau, X_{(k)}^{(-i)} < \tau \leq X_{(k+1)}^{(-i)}\},
\end{aligned}
$$

which is the required inequality in Theorem 2. $\square$

PROOF OF LEMMA 2.   The FNP is given by

$$
\begin{aligned}
N(t; J_0, J_1) &= \sum_{i \in J_1} \sum_{j=1}^{n} \frac{1}{j} I\{A = j, X_i < t\} \\
&= \sum_{i \in J_1} \sum_{j=1}^{n} \frac{1}{j} I\{X_{(j)} < t \le X_{(j+1)}, X_i < t\} \\
&= \sum_{i \in J_1} \sum_{j=1}^{n} \frac{1}{j} [I\{X_{(j-1)}^{(-i)} < t, X_i < t\} - I\{X_{(j)}^{(-i)} < t, X_i < t\}] \\
&= \sum_{i \in J_1} I\{X_i < t\} - \sum_{i \in J_1} \sum_{j=1}^{n-1} \frac{I\{X_{(j)}^{(-i)} < t, X_i < t\}}{j(j+1)}.
\end{aligned}
$$
(7.6)

Taking the expectation of (7.6), we get the first expression of the FNR. The second expression follows from the fact that

$$
\begin{aligned}
N(t; J_0, J_1) &= I\{A > 0\} - I\{UI(A > 0)/A\} \\
&= I\{X_{(1)} < t\} - \sum_{i \in J_0} \sum_{j=1}^{n} \frac{1}{j} I\{A = j, X_i < t\}.
\end{aligned}
$$
(7.7)        □

PROOF OF THEOREM 4.   We have

$$
\begin{aligned}
\mathrm{FNR}_\theta^{(2)}&(t_\tau \le \tau; J_0, J_1) \\
&= \sum_{k=1}^{n} E_\theta\{[N(t_\tau(k); J_0, J_1)|X_{(k)} < \tau \le X_{(k+1)}] \\
&\qquad\qquad\qquad \times I\{X_{(k)} < \tau \le X_{(k+1)}\}\} \\
&\le \sum_{i \in J_1} \sum_{k=1}^{n} \left\{ \frac{1}{k} \left[ 1 - \left( 1 - \frac{F_0(t_\tau(k))}{F_0(\tau)} \right)^k \right] \right. \\
&\qquad\qquad\qquad \left. \times P_\theta\{X_i < \tau, X_{(k-1)}^{(-i)} < \tau \le X_{(k)}^{(-i)}\} \right\}.
\end{aligned}
$$
(7.8)

The inequality in (7.8) follows from Theorem 3, noting that the conditional expectation of $N(t_\tau(k); J_0, J_1)$, given $\{X_{(k)} < \tau \le X_{(k+1)}\}$, is the FNR of the single-step procedure based on independent $Z_1, \ldots, Z_k$ with $Z_i \sim f_{\theta_i}(x)I(x < \tau)/F_{\theta_i}(\tau)$. The required inequality in Theorem 4 then follows from (7.8) because $F_{\theta_i}(\tau)$ is decreasing in $\theta_i$ for $i \in J_1$.   □

PROOF OF THEOREM 5. Since $V = \sum_{i=1}^{n} I(X_i \geq t) I(H_i = 0)$, we first note from Lemma 1 that the FDR under the mixture model is given by

$$
\begin{aligned}
\text{FDR}(t, n) &= \sum_{i=1}^{n} E_{\mathbf{H}} \bigg[ P\{X_i \geq t | H_i = 0\} \\
&\qquad - \sum_{j=1}^{n-1} \frac{P\{X_{(j)}^{(-i)} \geq t, X_i \geq t | \mathbf{H} \text{ with } H_i = 0\}}{(n-j)(n-j+1)} \bigg] \\
(7.9) \qquad &= \sum_{i=1}^{n} \bigg[ P\{X_i \geq t, H_i = 0\} - \sum_{j=1}^{n-1} \frac{P\{X_{(j)}^{(-i)} \geq t, X_i \geq t, H_i = 0\}}{(n-j)(n-j+1)} \bigg] \\
&= \sum_{i=1}^{n} \bigg[ P\{X_i \geq t, H_i = 0\} \\
&\qquad \times \bigg\{ 1 - \sum_{j=1}^{n-1} \frac{P\{X_{(j)}^{(-i)} \geq t | X_i \geq t, H_i = 0\}}{(n-j)(n-j+1)} \bigg\} \bigg].
\end{aligned}
$$

We now prove that

$$
(7.10) \qquad P\{X_{(j)}^{(-i)} \geq t | X_i \geq t, H_i = 0\} \geq P\{X_{(j)}^{(-i)} \geq t | X_i \geq t\}
$$

under the assumed positive dependence condition of the density $f$ of $\mathbf{X}$.

Let $\psi(X_i) = P\{X_{(j)}^{(-i)} \geq t | X_i, \theta_i = 0\}$. Then the conditional probability $P\{X_{(j)}^{(-i)} \geq t | X_i \geq t, \theta_i\}$ can be written as

$$
(7.11) \qquad \frac{E\{\psi(X_i) I(X_i \geq t - \theta_i)\}}{E\{I(X_i \geq t - \theta_i)\}},
$$

with the expectations taken with respect to $X_i$ under $\theta_i = 0$. Note that $\psi(x)$ is an increasing function of $x$ under the assumed positive dependence condition of $f$. Also, $I(x \geq t)$ is a totally positive of order two (TP$_2$) function of $(x, t)$ (see, e.g., [11]). Therefore, the ratio

$$
(7.12) \qquad \frac{E\{\psi(X_i) I(X_i \geq t)\}}{E\{I(X_i \geq t)\}}
$$

is increasing in $t$, because it is the expectation of an increasing function of a random variable whose distribution is stochastically increasing in $t$. This proves that $P\{X_{(j)}^{(-i)} \geq t | X_i \geq t, H_i = 0\} \geq P\{X_{(j)}^{(-i)} \geq t | X_i \geq t, H_i = 1\}$, implying that the probability $P\{X_{(j)}^{(-i)} \geq t | X_i \geq t\}$, being a convex combination of $P\{X_{(j)}^{(-i)} \geq t | X_i \geq t, H_i = 0\}$ and $P\{X_{(j)}^{(-i)} \geq t | X_i \geq t, H_i = 1\}$, is less than or equal to $P\{X_{(j)}^{(-i)} \geq t | X_i \geq t, H_i = 0\}$. Thus the required inequality (7.10) follows.

Applying (7.10) to (7.9), we get the inequality (6.2) to be proved in the theorem. The fact that

$$(7.13) \qquad \sum_{i=1}^{n} \delta_i = P\left\{ \max_{1 \leq i \leq n} X_i \geq t \right\} = P\{R > 0\}$$

follows from [13]. Furthermore, it is clear that the equality in (6.2) holds under independence of $(X_i, H_i)$. Thus, the theorem is proved. □

**8. Concluding remarks.** We have obtained in this article some theoretical results that extend previous work done under the assumption of independent tests. Two of these set the stage for developing our idea to modify the FDR- and FNR-controlling Bonferroni and Šidák procedures and obtaining wider families of FDR- and FNR-controlling procedures. We developed this idea by extending inequalities for FDR and FNR under independence from single-step to two-step procedures. In the case of the Bonferroni procedures, it is somewhat similar to what Storey, Taylor and Siegmund [19] used to modify the FDR-controlling BH procedure (which is, of course, a stepwise procedure) under independence. In the case of Šidák procedures, however, it is stronger in that we consider modifying less conservative procedures. It is important to point out that modifying the Šidák procedure by simply finding $t$ that controls $(\hat{n}_0/n)\{1 - F_0^n(t)\}$ (the estimated maximum FDR, which is basically the idea in modifying the FDR-controlling Bonferroni procedure), does not seem to provide much improvement to the Šidák procedure. The same is true for the FNR-controlling Šidák procedure. This is what we have noticed based on additional simulations not reported here. Also, as is seen from Table 2, we need to be cautious using the present modifications when there is too much dependence in the tests; they may become anticonservative. Procedures that control FDR are different from those that control FNR. It will be interesting to see if procedures that control both FDR and FNR can be developed using the results discussed in this paper.

**Acknowledgments.** The author thanks an Associate Editor and two referees for valuable comments, which have resulted in an improved paper, and Tianhui Zhou for her help with the numerical calculations.

FOX SCHOOL OF BUSINESS AND MANAGEMENT
TEMPLE UNIVERSITY
SPEAKMAN 319
1810 NORTH 13TH STREET
PHILADELPHIA, PENNSYLVANIA 19122-6083
USA
E-MAIL: sanat@temple.edu